\documentclass[draft]{amsart}

\title[Truncations of $p$-divisible groups]{Minimal truncations of supersingular $p$-divisible groups}
\author{Marc-Hubert Nicole and Adrian Vasiu}
\date{\today}

\usepackage{amssymb}
\usepackage{graphicx, epsfig}
\usepackage{
latexsym, amsfonts, amscd, amsmath}

\newtheorem{thm}{Theorem}[section]
\newtheorem{prop}[thm]{Proposition}
\newtheorem{lem}[thm]{Lemma}
\newtheorem{cor}[thm]{Corollary}

\newtheorem{con}[thm]{Conjecture}

\newtheorem{rmk}[thm]{Remark}

\newtheorem{exa}[thm]{Example}
\newtheorem{sch}[thm]{Scholium}

\input amssym.def
\newcommand{\calA}{{\mathcal{A}}}

\newcommand{\calG}{{\mathcal{G}}}

\newcommand{\calI}{{\mathcal{I}}}

\newcommand{\calV}{{\mathcal{V}}}

\newcommand{\End}{{\operatorname{End}}}

\newcommand{\Hom}{{\operatorname{Hom}}}

\newcommand{\Ker}{{\operatorname{Ker}}}

\newcommand{\Spec}{{\operatorname{Spec }}}

\def\db#1{{\fam\msbfam\relax#1}}

 \def\dbF{{\db F}}

 \def\dbN{{\db N}}
 
\def\dbQ{{\db Q}}

 \def\dbZ{{\db Z}}

\def\vph{{\varphi}}

\def\Ker{\text{Ker}}

\def\Hom{\text{Hom}}
\def\End{\text{End}}
\def\Spec{\text{Spec}}

\def\Lie{\text{Lie}}

\def\leaderfill{\leaders\hbox to 1em
     {\hss.\hss}\hfill}

\def\finishproclaim{\par\rm
     \ifdim\lastskip\medskipamount\removelastskip
     \penalty55\medskip\fi}
\def\endproof{$\hfill \square$}
\def\proof{\par\noindent {\it Proof:}\enspace}

\def\Ref[#1]{\par\hang\indent\llap{\hbox to\parindent
     {[#1]\hfil\enspace}}\ignorespaces}
\def\Item#1{\par\smallskip\hang\indent\llap{\hbox to\parindent
     {#1\hfill$\,\,$}}\ignorespaces}
\def\ItemItem#1{\par\indent\hangindent2\parindent
     \hbox to \parindent{#1\hfill\enspace}\ignorespaces}

\def\Le{{\mathchoice{\,{\scriptstyle\le}\,}
  {\,{\scriptstyle\le}\,}
  {\,{\scriptscriptstyle\le}\,}{\,{\scriptscriptstyle\le}\,}}}
\def\Ge{{\mathchoice{\,{\scriptstyle\ge}\,}
  {\,{\scriptstyle\ge}\,}
  {\,{\scriptscriptstyle\ge}\,}{\,{\scriptscriptstyle\ge}\,}}}

\def\arrowsim{\,\smash{\mathop{\to}\limits^{\lower1.5pt
  \hbox{$\scriptstyle\sim$}}}\,}

\def\doublemaprights#1#2#3#4{\raise3pt\hbox{$\mathop{\,\,\hbox to
#1pt{\rightarrowfill}\kern-30pt\lower3.95pt\hbox to
     #2pt{\rightarrowfill}\,\,}\limits_{#3}^{#4}$}}

\def\rightcapdownarrow{\raise9pt\hbox{$\ssize\cap$}\kern-7.75pt
     \Big\downarrow}

\def\rcapmapdown#1{\rightcapdownarrow\kern-1.0pt\vcenter{
     \hbox{$\scriptstyle#1$}}}

\def\rmapdown#1{\Big\downarrow\kern-1.0pt\vcenter{
     \hbox{$\scriptstyle#1$}}}
\def\rightsubsetarrow#1{{\ssize\subset}\kern-4.5pt\lower2.85pt
     \hbox to #1pt{\rightarrowfill}}
\def\longtwoheadedrightarrow#1{\raise2.2pt\hbox to #1pt{\hrulefill}
     \!\!\!\twoheadrightarrow}

\def\Hom{\operatorname{\hbox{Hom}}}

\usepackage[usenames]{color}
\definecolor{Indigo}{rgb}{0.2,0.1,0.7}
\definecolor{Violet}{rgb}{0.5,0.1,0.7}
\definecolor{White}{rgb}{1,1,1}
\definecolor{Green}{rgb}{0.1,0.9,0.2}

\begin{document}

\maketitle

\medskip\noindent
{\bf ABSTRACT.} Let $k$ be an algebraically closed field of
characteristic $p>0$. Let $H$ be a supersingular $p$-divisible group
over $k$ of height $2d$. We show that $H$ is uniquely determined
up to isomorphism by its truncation of level $d$ (i.e., by
$H[p^d]$). This proves Traverso's truncation conjecture for
supersingular $p$-divisible groups. If $H$ has a principal
quasi-polarization $\lambda$, we show that
$(H,\lambda)$ is also uniquely determined up to isomorphism
by its principally quasi-polarized truncated Barsotti--Tate group
of level $d$ (i.e., by $(H[p^d],\lambda[p^d])$).


\section{Introduction}

Let $p\in\dbN$ be a prime. Let $k$ be an algebraically closed
field of characteristic $p$. Let $c,d\in\dbN$. Let $H$ be a
$p$-divisible group over $k$ of codimension $c$ and dimension $d$;
thus the height of $H$ is $c+d$. Let $n\in\dbN$ be the smallest
number such that $H$ is uniquely determined up to isomorphism by
$H[p^n]$ (i.e., if $H_1$ is a $p$-divisible group over $k$ such
that $H_1[p^n]$ is isomorphic to $H[p^n]$, then $H_1$ is
isomorphic to $H$). Thus $H[p^n]$ is the {\it minimal truncation}
of $H$ which determines $H$. It is known that the number $n$
admits upper bounds that depend only on $c$ and $d$ (see [Ma],
[Tr1, Thm. 3], [Tr2, Thm. 1], [Va1, Cor. 1.3], and [Oo, Cor. 1.7]). For instance, Traverso
proved that $n\Le cd+1$ (see [Tr1, Thm. 3]). Traverso's work on
Grothendieck's specialization conjecture led him to speculate that
much more is true (cf. [Tr3, \S40, Conj. 4]):

\begin{con} \label{traversocon} We have $n\Le\min\{c,d\}$.
\end{con}

\indent We suppose for the remainder of the paper that the
codimension and the dimension of $H$ are equal. We recall that $H$
is called {\it supersingular} if all the slopes of its Newton
polygon are ${1\over 2}$. We prove the Conjecture for the
supersingular case:

\begin{thm} \label{traversosupersingular} Suppose $H$ is a
supersingular $p$-divisible group over $k$ of height $2d$. Then
$n\Le d$ i.e., $H$ is uniquely determined up to isomorphism by
$H[p^d]$.
\end{thm}

\noindent Theorem \ref{traversosupersingular} strengthens
Traverso's and Vasiu's results (see [Tr1, Thm. 3] and [Va1, Prop. 4.1.1])
which worked with $H[p^{d^2+1}]$ and $H[p^{d^2}]$ (respectively).
Theorem 1.2 was originally claimed in [Ni, \S 1.4]. It turns out
that [Ni, Lem. 1.4.4, Cor. 1.4.7] are incorrect as stated. The
proof in the present paper uses elementary methods of
$\sigma$-linear algebra to avoid those issues completely.

For the sake of completeness, we also prove the following
principally quasi-polarized variant of Theorem
\ref{traversosupersingular}.

\smallskip
\begin{thm} \label{traversopp} Suppose $H$ is a supersingular $p$-divisible group over $k$ of height $2d$ which has a principal quasi-polarization
$\lambda$. Then $(H,\lambda)$ is uniquely determined up to isomorphism by $(H[p^d], \lambda[p^d])$ (i.e., by its principally quasi-polarized truncated Barsotti--Tate group of level $d$).
\end{thm}

\smallskip
\noindent Theorem \ref{traversopp} refines and extends [Va1, Prop.
4.1.1]. Theorems \ref{traversosupersingular} and \ref{traversopp}
 are optimal (i.e., they do not hold if $[p^d]$ gets replaced
with $[p^{d-1}]$; see Example \ref{effectivebound}). Theorems 1.2 and 1.3 and Corollary 3.2 represent progress towards:

\medskip
(i) the classification of (principally quasi-polarized) supersingular $p$-divisible groups, and

\smallskip
(ii) the understanding of the ultimate stratifications introduced
in [Va1, Basic Thm. 5.3.1 and Subsubsection 5.3.2] and of the
level $m$ stratifications introduced in [Va2, Subsections 4.2 and
4.4].

\medskip\indent
In Section 2, we introduce notations and basic invariants
which pertain to supersingular $p$-divisible groups and which
allow us to get a more precise form of Theorem
\ref{traversosupersingular} (see Corollary 3.2). In Sections 3 and
4, we prove Theorems \ref{traversosupersingular} and
\ref{traversopp} (respectively). Note that the proof of
Theorem \ref{traversosupersingular} is self-contained.

\section{Basic invariants of supersingular Dieudonn\'e modules}

Let $W(k)$ be the ring of Witt vectors with coefficients in $k$.
For $s\in\dbN$, let $W_s(k):=W(k)/(p^s)$. Let $\sigma$ be the
Frobenius automorphism of $W(k)$ (or
$W_s(k)$) induced from $k$.

\begin{prop} \label{wellknown} Let $\calG$ be a smooth group scheme over
$\Spec(\dbZ_p)$ such that its special fibre $\calG_{\dbF_p}$ is a
connected, affine scheme. Let $\sigma$ act naturally on
$\calG(W(k))$. Then we have
$\calG(W(k))=\{g_0^{-1}\sigma(g_0)|g_0\in \calG(W(k))\}$.
\end{prop}

\medskip
\proof Let $g\in\calG(W(k))$. By induction on $s\in\dbN$, we check
that there exists an element $g_s\in\calG(W(k))$ such that the
following two properties hold:

\medskip
(i) we have $g_sg\sigma(g_s)^{-1}\in\Ker(\calG(W(k))\to\calG(W_s(k)))$, and

\smallskip
(ii) for $s\ge 2$, the images of $g_s$ and $g_{s-1}$ in
$\calG(W_{s-1}(k))$ coincide.

\medskip
Let $\sigma$ act naturally on
$\calG(W_s(k))$. As $\calG_{\dbF_p}$ is affine and connected,
there exists $\bar g_1\in \calG(k)$ such that $\bar
g_1^{-1}\sigma(\bar g_1) $ is the reduction mod $p$ of $g$ (cf.
Lang's theorem in [Bo, Ch. V, Cor. 16.4]). If $g_1\in\calG(W(k))$
lifts $\bar g_1$, then we have
$g_1g\sigma(g_1)^{-1}\in\Ker(\calG(W(k))\to\calG(k))$. Thus the
basis of the induction holds. The passage from $s$ to $s+1$ goes
as follows. As $\calG$ is smooth, the group
$\Ker(\calG(W_{s+1}(k))\to\calG(W_s(k)))$ is the group of
$k$-valued points of the vector group $\calV$ over $\dbF_p$
defined by $\Lie(\calG_{\dbF_p})$. From Lang's theorem applied to
$\calV$, we get that there exists $\bar
g^{\prime}_{s+1}\in\Ker(\calG(W_{s+1}(k))\to\calG(W_s(k)))$ such
that $(\bar g^{\prime}_{s+1})^{-1}\sigma(\bar g^\prime_{s+1})$ is
the image of $g_sg\sigma(g_s)^{-1}$ in
$\Ker(\calG(W_{s+1}(k))\to\calG(W_s(k)))$. Let
$g^\prime_{s+1}\in\calG(W(k))$ be an element that lifts $\bar
g^\prime_{s+1}$. If $g_{s+1}:=g^\prime_{s+1}g_s\in\calG(W(k))$,
then we have $g_{s+1}g\sigma(g_{s+1})^{-1}$
$\in\Ker(\calG(W(k))\to\calG(W_{s+1}(k)))$. Moreover, as $\bar
g^{\prime}_{s+1}\in\Ker(\calG(W_{s+1}(k))\to\calG(W_s(k)))$, the
images of $g_{s+1}$ and $g_s$ in $\calG(W_s(k))$ coincide. This
ends the induction.

Due to the property (ii), the $p$-adic limit of the sequence $(g_s)_{s\in\dbN}$
is an element $g_{\infty}$ in $\calG(W(k))$. Due to the property (i), the
element $g_{\infty}g\sigma(g_{\infty})^{-1}$ is the identity. Thus
$g_{\infty}^{-1}\sigma(g_{\infty})=g$. This implies that
$\calG(W(k))=\{g_0^{-1}\sigma(g_0)|g_0\in \calG(W(k))\}$.\endproof

\subsection{} Let $H$ be a supersingular $p$-divisible group over $k$ of height
$2d$, for $d \in \dbN$. Let $(M,\phi)$ be the (contravariant)
Dieudonn\'e module of $H$. Thus $M$ is a free $W(k)$-module of
rank $2d$ and $\phi:M\to M$ is a $\sigma$-linear endomorphism such
that we have $pM\subseteq \phi(M)$. Let $\vartheta:M\to M$ be the
Verschiebung map of $(M,\phi)$; we have
$\vartheta\phi=\phi\vartheta=p1_M$. We denote also by $\phi$ the
$\sigma$-linear automorphism of $\End(M[{1\over p}])$ that takes
$e\in \End(M[{1\over p}])$ to $\phi(e):=\phi\circ
e\circ\phi^{-1}\in \End(M[{1\over p}])$. Let $A:=\{e\in \End(M)|\phi(e)=e\}$ be the
$\db{Z}_p$-algebra of endomorphisms of $(M,\phi)$. Let $\calA$ be
the smooth, affine group scheme over $\Spec(\dbZ_p)$ of invertible
elements of the $\dbZ_p$-algebra $A$. Let $O$ be the $W(k)$-span
of $A$. As all slopes of $(M,\phi)$ are ${1\over 2}$, all slopes
of $(\End(M[{1\over p}]),\phi)$ are $0$. This implies that $O$ is
a $W(k)$-subalgebra of $\End(M)$ such that the quotient
$W(k)$-module $\End(M)/O$ is torsion. Let $E$ be the (unique up to
isomorphism) supersingular $p$-divisible group over $k$ of height
$2$. Let $(N,\vph)$ and $(N^d,\vph)$ be the Dieudonn\'e modules of
$E$ and $E^d$ (respectively).

Let $n\in\dbN$ be as in Section 1.
We list two additional basic invariants of $H$:

\medskip

$\bullet$ Let $m\in\dbN$ be the smallest number such that
$p^m\End(M)\subseteq O\subseteq\End(M)$.

\smallskip
$\bullet$  Let $q\in\dbN\cup\{0\}$ be the smallest number such
that there exists a monomorphism $j:(N^d,\vph)\hookrightarrow
(M,\phi)$ with the property that $\phi^q(M)\subseteq j(N^d)$.

\medskip\noindent
Thus $m$ is the
Fontaine--Dieudonn\'e torsion of $(\End(M),\phi)$ defined in [Va1, 2.2.2 (b)].
We have $q=0$ if and only if there exists an isomorphism $(N^d,\vph)\arrowsim (M,\phi)$.

\smallskip
\begin{thm} \label{supersingularinvariants}
Let $(M,\phi)$ be a supersingular Dieudonn\'e module over $k$, and
let $H$ be the corresponding $p$-divisible group. Let $n$, $m$,
and $O$ be as above.

\medskip
\noindent {\bf (a)} Let $t\in\dbN$ be the smallest number such
that for each element $g\in\pmb{GL}_M(W(k))$ congruent mod $p^t$
to $1_M$, there exists an isomorphism between the two Dieudonn\'e
modules $(M,g\phi)$ and $(M,\phi)$. Then $t=n$ i.e., $t$ is the
smallest number such that $H$ is uniquely determined up to
isomorphism by $H[p^t]$.

\smallskip
\noindent {\bf (b)} Let $g\in\pmb{GL}_M(W(k))\cap O$. Then $(M,g\phi)$ and $(M,\phi)$ are isomorphic.

\smallskip
\noindent {\bf (c)} If $p^m\End(M)\subseteq O\subseteq\End(M)$
then $H$ is determined up to isomorphism by $H[p^m]$ i.e., we have
an inequality $n\le m$.
\end{thm}

\medskip
\proof We first prove (a). This is a special case of [Va1, Lemma 3.2.2]
for the group $G= \pmb{GL}_M$, but for the sake of completeness we
include a self-contained proof which works for all $p$-divisible
groups over $k$. Let us first show that $t\Le n$. Let $g\in
\pmb{GL}_M(W(k))$ be congruent mod $p^n$ to $1_M$. Let $H_g$ be
the $p$-divisible group over $k$ whose Dieudonn\'e module is
$(M,g\phi)$. Then $H_g[p^n]=H[p^n]$ and thus $H_g$ and $H$ are
isomorphic i.e., $(M,\phi)$ and $(M,g\phi)$ are isomorphic. Thus
$t\Le n$.

\indent Second, we show that $n\Le t$. Let $H_t$ be a
$p$-divisible group over $k$ such that $H_t[p^t]$ and $H[p^t]$ are
isomorphic. Let $g\in \pmb{GL}_M(W(k))$ be such that the
Dieudonn\'e module of $H_t$ is isomorphic to $(M,g\phi)$. As
$H_t[p^t]$ and $H[p^t]$ are isomorphic, we can assume that
$(M,g\phi,\vartheta g^{-1})$ mod $p^t$  is $(M,\phi,\vartheta)$
mod $p^t$. This implies that $g$ fixes $\phi(M)/p^tM$ and
$M/p^{t-1}\phi(M)$. Since $g$ fixes
$pM/p^tM\subseteq\phi(M)/p^tM$, there exists $u\in\End(M)$ such
that $g=1_M+p^{t-1}u$. As $g$ fixes $\phi(M)/p^tM$ and
$M/p^{t-1}\phi(M)$, we get that $u$ mod $p$ annihilates
$\phi(M)/pM$ and $M/\phi(M)$. Thus $u(\phi(M))\subseteq pM$ and
$u(M)\subseteq\phi(M)$. This implies that $u^2\in p\End(M)$, that
$(\phi^{-1}u\phi)(M)\subseteq \phi^{-1}(pM)=\vartheta(M)$, and
that $(\phi^{-1}u\phi)(\vartheta(M))=\phi^{-1}(u(pM))\subseteq
pM$. Let $v:=\phi^{-1}u\phi\in\End(M)$; we have $u=\phi(v)$ and
$v$ mod $p$ fixes $\vartheta(M)/pM$ and $M/\vartheta(M)$. As
$\vartheta(M)/pM$ is the kernel of $\phi$ mod $p$, it is easy to
see that we can write $v=pv_1+v_2$, where $v_1,v_2\in\End(M)$ and
$\phi(v_2)\in p\End(M)$.

If $g^\prime\in\Ker(\pmb{GL}_M(W(k))\to\pmb{GL}_M(W_t(k)))$ and if
$(M,g^\prime g\phi)$ is isomorphic to $(M,\phi)$, then $(M,g\phi)$
is isomorphic to $(M,g^{\prime\prime}\phi)$ for some
$g^{\prime\prime}\in\Ker(\pmb{GL}_M(W(k))\to\pmb{GL}_M(W_t(k)))$;
thus $(M,g\phi)$ is also isomorphic to $(M,\phi)$ (cf. the
definition of $t$). Thus to show that $(M,g\phi)$ and $(M,\phi)$
(i.e., that $H_t$ and $H$) are isomorphic, we can replace $g$ by
any element of $\pmb{GL}_M(W(k))$ congruent mod $p^t$ to $g$. In
other words, we can replace $u$ by any element of $u+p\End(M)$. By
replacing $u$ with $u-\phi(v_2)$ and $v$ with $v_1=v-v_2$, we can
assume that $v=pv_1\in p\End(M)$. We define
$g_1:=(1_M-p^tv_1)^{-1}\in\Ker(\pmb{GL}_M(W(k))\to\pmb{GL}_M(W_t(k)))$
and $g_2:=g_1g\phi(g_1^{-1})$. We have
$g_2=g_1g\phi(1_M-p^tv_1)=g_1(1_M+p^{t-1}u)(1_M-p^{t-1}u)=g_1(1_M-p^{2t-2}u^2).$
As $u^2\in p\End(M)$ and $t\ge 1$, we have $p^{2t-2}u^2\in
p^t\End(M)$. Thus $g_2$ is congruent mod $p^t$ to $1_M$. From the
definition of $t$ we get
 that $(M,g_2\phi)$ and $(M,\phi)$ are isomorphic. As $g_2\phi=g_1g\phi g^{-1}_1$, we conclude that $(M,g\phi)$ and $(M,\phi)$ are isomorphic. Thus $H_t$ and $H$ are isomorphic. This implies that $n\le t$. Thus $n=t$ and therefore (a) holds.

\indent Part (b) is a particular case of [Va1, proof of Cor.
3.3.4], but we provide here a simpler argument which works for all
isoclinic $p$-divisible groups. The inverse in $\pmb{GL}_M(W(k))$
of the element $g\in O$ is a polynomial in $g$ with coefficients
in $W(k)$ (cf. the Cayley--Hamilton theorem) and thus it
belongs to $O$. Thus $g$ has an inverse in $O$ and therefore
$g\in\calA(W(k))$. Each invertible element of $O$ is also an
invertible element of $\End(M)$ and therefore we have
$\calA(W(k))\leqslant \pmb{GL}_M(W(k))$. The automorphism $\sigma$
acts naturally on $\calA(W(k))$. \noindent As $\calA$ is an open
subscheme of the vector group scheme over $\Spec(\dbZ_p)$ defined
by $A$, its fibres are connected. Thus there exists
$g_0\in\calA(W(k))$ such that $g_0^{-1}\sigma(g_0)=g$, cf.
Proposition \ref{wellknown}. We have $g_0g\sigma(g_0)^{-1}=1_M$.
As $\sigma(g_0)=\phi(g_0)$, we have $g_0g\phi g_0^{-1}=\phi$. Thus
$g_0$ is an isomorphism between $(M,g\phi)$ and $(M,\phi)$. Thus
(b) holds.

Based on (a), to prove (c) it suffices to show that for each
element $g\in \pmb{GL}_M(W(k))$ congruent mod $p^m$ to $1_M$, the
Dieudonn\'e modules $(M,g\phi)$ and $(M,\phi)$ are isomorphic. As
$g-1_M$ $\in p^m\End(M)\subseteq O$, we have $g\in O$. Thus
$(M,g\phi)$ and $(M,\phi)$ are isomorphic, cf. (b). Thus (c)
holds.\endproof

\begin{sch} \label{scholium1} \emph{Let $\{x,y\}$ be a $W(k)$-basis for $N$ such that
$\vph(x)=y$ and $\vph(y)=px$. Thus $\{px,y\}$ is a $W(k)$-basis
for $\vph(N)$ and we have $\phi(px)=py$ and
$\phi^2(px)=\phi(py)=p^2x$. The image of the map $\vph^2-p1_N:N\to
N$ is $pN$. Let $N^*:=\Hom(N,W(k))$. Let $\{x^*,y^*\}$ be the
$W(k)$-basis for $N^*$ which is the dual of $\{x,y\}$. Thus
$\{x\otimes x^*,y\otimes y^*,x\otimes y^*,y\otimes x^*\}$ is a
$W(k)$-basis for $\End(N)=N\otimes_{W(k)} N^*$. The
$\sigma$-linear automorphism $\vph$ of $\End(N[{1\over p}])$
permutes $x\otimes x^*$ and $y\otimes y^*$ as well as $px\otimes
y^*$ and $y\otimes x^*$. Thus $\{x\otimes x^*,y\otimes
y^*,px\otimes y^*,y\otimes x^*\}$ is a $W(k)$-basis for the
$W(k)$-span $O_1$ of endomorphisms of $(N,\vph)$. We have
inclusions $p\End(N)\varsubsetneqq O_1\varsubsetneqq \End(N)$.}

\indent \emph{Let $O_d$ be the $W(k)$-span of the
$\db{Z}_p$-algebra of endomorphisms of $(N^d,\phi)$. The inclusion
$O_d\subseteq\End(N^d)$ can be identified with the inclusion of
matrix $W(k)$-algebras $M_d(O_1)\subseteq M_d(\End(N))$. Thus we
have $p\End(N^d)\varsubsetneqq  O_d\varsubsetneqq \End(N^d)$.
\noindent If $H \cong E^d$, we thus retrieve the well-known result
that $H$ is determined by its $p$-kernel, since $n \leq m = 1$.}
\end{sch}

\begin{lem} \label{manin} Let $d$ be a natural number
and let $(N^d, \varphi)$ be as above. Let $(M,\phi)$ be a
supersingular Dieudonn\'e module of rank $2d$. Then there exists a
monomorphism $j: (N^d,\varphi) \hookrightarrow (M,\phi)$ for which
$\phi^{d-1}(M)$ is contained in $j(N^d)$ i.e., $q \leq d-1$.
\end{lem}

\medskip
\proof We prove the Lemma by induction on $d\in\dbN$. If $d=1$,
then $H$ is isomorphic to $E$ and thus $q=0$. Suppose $d\Ge
2$. We consider a short exact sequence
$$0\to (N,\vph)\to (M,\phi)\to (M_1,\phi_1)\to 0$$
of supersingular Dieudonn\'e modules over $k$. As the rank of
$M_1$ is $2d-2$, by induction there exists a monomorphism
$j_1:(N^{d-1},\vph)\hookrightarrow (M_1,\phi_1)$ such that
$\phi_1^{d-2}(M_1)\subseteq j_1(N^{d-1})$. Let $M_2$ be the
inverse image of $\phi_1(j_1(N^{d-1}))$ in $M$.

\indent We have a short exact sequence
$$0\to (N,\vph)\to (M_2,\phi)\to (\phi_1(j_1(N^{d-1})),\phi_1)\to 0\leqno (1)$$
of supersingular Dieudonn\'e modules over $k$. We check that the
short exact sequence (1) splits. The Dieudonn\'e module
$(\phi_1(j_1(N^{d-1})),\phi_1)$  is a direct sum of supersingular
Dieudonn\'e modules of rank $2$ which have $W(k)$-bases $\{x,y\}$
with the properties that: (i)  $\phi_1(x)=py$ and
$\phi_1^2(x)=p\phi_1(y)=px$, and (ii) $x\in pj_1(N^{d-1})\subseteq
pM_1$ (see Scholium \ref{scholium1}). Thus to check that (1)
splits, it suffices to show that for each such $W(k)$-basis
$\{x,y\}$, there exists an element $x_2\in M_2$ such that it maps to $x$
and moreover, we have ${1\over p}\phi(x_2)\in M_2$ and
$\phi^2(x_2)=px_2$. Let $x_1\in pM$ be such that it maps to $x$,
cf. (ii). Let $y_1:=\phi^2(x_1)-px_1$; it is an element of $pN$.
Let $y_2\in N$ be such that $\vph^2(y_2)-py_2=-y_1$ (see Scholium \ref{scholium1}). Let $x_2:=x_1+y_2$; it is an element of
$M_2$ that maps to $x$ and we have $\phi^2(x_2)-px_2=y_1-y_1=0$. As $x_1\in pM$ and $y_2\in pN$, we have
${1\over p}\phi(x_2)={1\over p}\phi(x_1)+{1\over p}\vph(y_2)\in
M$. As ${1\over p}\phi(x_2)$ maps to $y$, we have ${1\over
p}\phi(x_2)\in M_2$. Thus the element $x_2$ exists. As
$\phi_1^{d-2}(M_1)\subseteq j_1(N^{d-1})$, we have
$\phi_1^{d-1}(M_1)\subseteq \phi_1(j_1(N^{d-1}))$. This implies
that $\phi^{d-1}(M)\subseteq M_2$. As the short exact sequence (1)
splits, there exists an isomorphism $j_2:(N^d,\vph)\arrowsim
(M_2,\phi)$. Its composite with the monomorphism
$(M_2,\phi)\hookrightarrow (M,\phi)$ is a monomorphism
$j:(N^d,\vph)\hookrightarrow (M,\phi)$ such that we have
$\phi^{d-1}(M)\subseteq j(N^d)=M_2$. Thus $q\Le d-1$. This ends
the induction.\endproof

\begin{rmk}
\emph{Lemma \ref{manin} follows also from [Ma, Thm. 3.7].}
\end{rmk}

\begin{rmk} \emph{\label{trivial0} The smallest number $\kappa\in\dbN\cup\{0\}$ such that there exists an isogeny $H\to E^d$ whose kernel is annihilated by $p^{\kappa}$, is $\lceil {q\over 2}\rceil$ (i.e., it is the smallest number such that $p^{\kappa}$ annihilates $N^d/\vph^q(N^d)$).
}
\end{rmk}

\begin{sch} \label{scholium2} \emph{For $i\in\dbN\cup\{0\}$, let
$f(i)$ be the biggest integer such that $M\subseteq
p^{f(i)}\phi^i(M)$. We have
$$O=\bigcap_{i=0}^{\infty}
\phi^i(\End(M))=\bigcap_{i=0}^{\infty}
\End(\phi^i(M))=\bigcap_{i=0}^{\infty} \End(p^{f(i)}\phi^i(M)).$$
Thus the Fontaine-Dieudonn\'e torsion $m\in\dbN$ is the smallest
number such that
$$M\subseteq \bigcup_{i \in \dbN} p^{f(i)}\phi^i(M)\subseteq
p^{-m}M.$$}

\end{sch}

\section{Proof of Theorem \ref{traversosupersingular}}

We recall that $d\in\dbN$, that $H$ is a supersingular $p$-divisible group over $k$ of height $2d$, that $(M,\phi)$ is the Dieudonn\'e module of $H$, and that we have introduced three invariants $n$, $m$, and $q$ of $H$ (see Subsection 2.1).

\begin{thm} \label{importantbound}
We have inclusions $p^{q+1}\End(M)\subseteq O\subseteq\End(M)$
i.e., $m \leq q+1$.
\end{thm}

\medskip
\proof We prove the Theorem by a step $2$ induction on $q\in\dbN
\cup\{0\}$. If $q=0$, then $H$ is isomorphic to $E^d$ and thus
$m=1=q+1$ (cf. Scholium \ref{scholium1}).

\indent Let $q=1$. Let $j:(N^d,\vph)\hookrightarrow
(M,\phi)$ be a monomorphism such that $\phi(M)\subseteq j(N^d)$.
We have $j(N^d)\subseteq M\subseteq \phi^{-1}(j(N^d))$. This
implies that we have a direct sum decomposition $j(N^d)=X\oplus
Y_1\oplus Y_2$ such that $M=X\oplus {1\over p}Y_1\oplus Y_2$,
$\phi(X)=Y_1\oplus Y_2$, and $\phi(Y_1\oplus Y_2)=pX$. Let
$i\in\dbN$. If $i$ is even, then $p^{{-i-2}\over
2}\phi^i(M)={1\over p}X\oplus {1\over p^2}Y_{1i}\oplus {1\over
p}Y_{2i}$, where $Y_{1i}:=p^{-{{i}\over 2}}\phi^i(Y_1)$ and
$Y_{2i}:=p^{-{{i}\over 2}}\phi^i(Y_2)$. As $Y=Y_{1i}\oplus
Y_{2i}$, we have $M\subseteq p^{{-i-2}\over 2}\phi^i(M)\subseteq
p^{-2}M$. If $i=2l+1$ is odd, then $p^{-l-i}\phi^i(M)={1\over
p}X_{1i}\oplus X_{2i}\oplus {1\over p}Y_1\oplus {1\over p}Y_2$,
where $X_{1i}:=p^{-l-1}\phi^i(Y_1)$ and
$X_{2i}:=p^{-l-1}\phi^i(Y_2)$. As $X=X_{1i}\oplus X_{i2}$, we have
$M\subseteq p^{-l-1}\phi^i(M)\subseteq p^{-1}M$. Regardless of
what $i\in\dbN$ is, we have $M\subseteq \cup_{i\in\dbN}
p^{f(i)}\phi^i(M)\subseteq p^{-2}M$ and thus $m\le 2=q+1$ (cf.
Scholium \ref{scholium2}).

\indent Suppose $q\ge 2$. Let $j:(N^d,\vph)\hookrightarrow
(M,\phi)$ be a monomorphism such that $\phi^q(M)\subseteq j(N^d)$.
Thus $\phi^{q-2}(M)\subseteq \phi^{-2}(j(N^d))={1\over p}j(N^d)$.
Let $\tilde M:={1\over p}j(N^d)+M$. Let $\tilde O$ be the
$W(k)$-subalgebra of $\End(\tilde M)$ generated by endomorphisms
of $(\tilde M,\phi)$. We have $\phi^{q-2}(\tilde
M)=\phi^{q-2}(j(N^d))+\phi^{q-2}(M)\subseteq j(N^d)+{1\over
p}j(N^d)\subseteq {1\over p}j(N^d)$. Let $\tilde
\j:(N^d,\phi)\hookrightarrow (\tilde M,\phi)$ be a monomorphism
whose image is ${1\over p}j(N^d)$. We have $\phi^{q-2}(\tilde
M)\subseteq \tilde \j(N^d)\subseteq\tilde M$. Thus by induction,
we have $p^{q-1}\End(\tilde M)\subseteq\tilde O$. As $M\subseteq
\tilde M\subseteq {1\over p}M$, we have
$p\End(M)\subseteq\End(\tilde M)\subseteq {1\over p}\End(M)$. This
implies that
$$p^{q+1}\End(M)\subseteq p^q\End(\tilde M)\subseteq p\tilde
O\subseteq p\End(\tilde M)\subseteq\End(M).$$
\noindent As $p\tilde O$ is
$W(k)$-generated by elements fixed by $\phi$ and as $p\tilde
O\subseteq\End(M)$, we have $p\tilde O\subseteq O$. Thus
$p^{q+1}\End(M)\subseteq p\tilde O\subseteq O$. This implies that
$m\Le q+1$. This ends the induction. \endproof

\medskip
\indent From Theorem \ref{supersingularinvariants} (c), Theorem
\ref{importantbound}, and Lemma \ref{manin} we get:

\begin{cor} \label{finaltouch} For each supersingular $p$-divisible group $H$ over $k$ of height $2d$, we have inequalities $n\leq m\leq q+1\leq d$.
\end{cor}

\medskip
\noindent This implies $n\leq d$ and ends the proof of Theorem
\ref{traversosupersingular}.

\begin{exa} \label{effectivebound} \emph{Let $d\Ge 2$. Suppose there exists a $W(k)$-basis
$\{e_1,\ldots,e_{2d}\}$ for $M$ such that for
$i\in\{1,\ldots,d\}$, we have $\phi(e_i)=e_{i+1}$ and for
$i\in\{d+1,\ldots,2d\}$, we have $\phi(e_i)=pe_{i+1}$ (here
$e_{2d+1}:=e_1$). We denote the corresponding $p$-divisible group
by $C_d$. Let $(M,\phi_1)$ be the Dieudonn\'e module with the
property that $\phi_1(e_i)=\phi(e_i)$ if $i\neq d+1$ and
$\phi_1(e_{d+1})=\phi_1^{d+1}(e_1)=pe_{d+2}+p^{d-1}e_2$. Let $H_1$
be the $p$-divisible group over $k$ whose Dieudonn\'e module is
$(M,\phi_1)$. We have $\phi_1^{2d}(e_1)=p^de_1+p^{d-1}e_{d+1}\in
p^{d-1}M\setminus p^dM$. But $\phi^{2d}(M)=p^dM$. From the last
two sentences, we get that $(M,\phi_1)$ is not isomorphic to
$(M,\phi)$ (i.e., $H_1$ is not isomorphic to $C_d$). It is easy to
see that $\phi_1$ and $\vartheta_1:=p\phi_1^{-1}$ are congruent
mod $p^{d-1}$ to $\phi$ and $\vartheta:=p\phi^{-1}$
(respectively). Thus $C_d[p^{d-1}]=H_1[p^{d-1}]$. From the last
two sentences, we get that $C_d$ is not determined by
$C_d[p^{d-1}]$. Thus $n\Ge d$. From this and Corollary
\ref{finaltouch}, we obtain the equalities $n=m=q+1=d$; thus the
inequalities of Corollary \ref{finaltouch} are best possible. }

\emph{We now discuss the polarized case. Let $\theta$ be an
invertible element of $W(k)$ such that we have
$\sigma^d(\theta)=-\theta$. Let $\psi:M\otimes_{W(k)} M\to W(k)$
be the perfect, alternating form on $M$ such that the following
two properties hold: (i) for $i,j\in\{1,\ldots,2d\}$ with
$|j-i|\neq d$, we have $\psi(e_i,e_j)=0$, and (ii) for
$i\in\{1,\ldots,d\}$ we have
$\psi(e_i,e_{i+d})=-\psi(e_{i+d},e_i)=\sigma^{i-1}(\theta)$. It is
easy to see that $\psi$ is a principal quasi-polarization of both
$(M,\phi)$ and $(M,\phi_1)$. Thus, if $\lambda$ is the principal
quasi-polarization of $C_d$ defined by $\psi$, then
$(C_d,\lambda)$ is not determined by
$(C_d[p^{d-1}],\lambda[p^{d-1}])$.}
\end{exa}

\begin{rmk} \emph{\label{trivial1} \label{trivial2} If $s\in\{2,\ldots,d\}$ and $H\cong C_s\times E^{d-s}$, then
$q=s-1$ (cf. Example \ref{effectivebound}). Thus $q$ can be any
number in the set $\{0,\ldots,d-1\}$. If $d=2\ell$ is even and
$H\cong C_2^{\ell}$, then $q=1$ and the $a$-number $\dim_k
(Hom(\pmb{\alpha}_p,H))$ of $H$ is $a=\ell$;
thus the difference $d-q-a=\ell-1$ can be any non-negative integer.}
\end{rmk}

\begin{rmk} \emph{\label{trivial3} Let $c',d'\in\dbN$ be relatively prime. Let $\ell\in\dbN$. Let $H'$ be a $p$-divisible group over $k$ of height $\ell(c'+d')$ and unique Newton polygon slope $\alpha:={d'\over {c'+d'}}$. If either $c'=1$ or $d'=1$, then the methods of this paper apply entirely to get an analogue of Corollary 3.2 for the slope $\alpha$ (and in particular, that $H'$ is uniquely determined up to isomorphism by $H'[p^{\ell\min\{c',d'\}}]$). Suppose $c',d'\Ge 2$ and $\ell=1$. The classical description of isogenies between such $p$-divisible groups $H'$ shows that the analogue of the invariant $q$ is an invariant $b$ which can be any number in the set $\{0,\ldots,(c'-1)(d'-1)\}$ (see [dJO, Subsections 5.8 and 5.32]). Moreover, the analogue of $\lceil{q\over 2}\rceil$ (see Remark \ref{trivial0}) is then $\lceil{b\over {c'+d'}}\rceil$. If $|d^\prime-c^\prime|\ge 3$, the inequality $\lceil{{2(c'-1)(d'-1)}\over {c'+d'}}\rceil+1>\min\{c',d'\}$ holds and therefore the mentioned description does not suffice to show that $H'$ is
uniquely determined up to isomorphism by $H'[p^{\min\{c',d'\}}]$.
The same applies if $|d^\prime-c^\prime|\in\{1,2\}$ and
$\ell>>0$.}
\end{rmk}

\section{Proof of Theorem \ref{traversopp}}

\subsection{} Let $H$ be a supersingular $p$-divisible group over $k$
of height $2d$ which has a principal quasi-polarization $\lambda$.
Let $(M,\phi)$, $A$, and $\calA$  be as in Subsection 2.1. Let $\psi$
be the perfect alternating form on $M$ induced by $\lambda$. Let
$\iota$ be the involution of $\End(M)$ defined by $\psi$: for
$x,y\in M$ and $e\in\End(M)$, we have an identity
$\psi(e(x),y)=\psi(x,\iota(e)(y))$. An element $e\in\End(M)$
annihilates $\psi$ (i.e., for all $x,y\in M$ we have
$\psi(e(x),y)+\psi(x,e(y))=0$) if and only if $\iota(e)=-e$.

Let
$G:=\pmb{Sp}(M,\psi)$; it is a reductive, closed subgroup scheme of
$\pmb{GL}_M$ whose Lie algebra $\Lie(G)$ is
$\{e\in\End(M)|\iota(e)=-e\}$. Moreover, for
an element $g\in\pmb{GL}_M(W(k))$, we have $g\in G(W(k))$ if and only if
$\iota(g)g=1_M$.

For $x,y\in M$, we have
$\psi(\phi(x),\phi(y))=p\sigma(\psi(x,y))$. This implies that
$\iota(A)=A$. It also implies that $\phi$ normalizes the Lie
subalgebra $\Lie(G)[{1\over p}]$ of $\End(M[{1\over p}])$. Thus
the triple $(M,\phi,G)$ is a latticed $F$-isocrystal with a group
over $k$ as defined in [Va1, 1.1 (a)]. As $\iota(A)=A$, the
involution $\iota$ acts naturally on all points of $\calA$ with
values in $\dbZ_p$-algebras. Let $\calI_{\dbQ_p}$ be the closed
subgroup of $\calA_{\dbQ_p}$ with the property that for each
$\dbQ_p$-algebra $R$, we have
$\calI_{\dbQ_p}(R)=\{g\in\calA(R)|\iota(g)g=1_{M\otimes_{\dbQ_p}
R}\}$. Let $\calI$ be the Zariski closure of $\calI_{\dbQ_p}$ in
$\calA$; it is a flat, closed subgroup scheme of $\calA$ whose
generic fibre is $\calI_{\dbQ_p}$.

\begin{lem} \label{p3lemma} Suppose that $p>2$. Then $\calI$ is a smooth group scheme over $\Spec(\dbZ_p)$.
\end{lem}

\medskip
\proof Let $B(k)$ be the field of fractions of $W(k)$. As
$G_{B(k)}=\calI_{B(k)}$, the group $\calI_{\dbQ_p}$ is connected.
Let $A^-:=\{e\in A|\iota(e)=-e\}$ and $A^+:=\{e\in
A|\iota(e)=e\}$. As $p>2$ and $\iota^2$ is the identity
automorphism of $A$, we have a direct sum decomposition
$A=A^-\oplus A^+$ of $\dbZ_p$-modules. The Lie algebra
$\Lie(\calI_{\dbF_p})$ is included in $A^-/pA^-$ and thus its
dimension is at most equal to the dimension of
$A^-\otimes_{\dbZ_p} B(k)=\Lie(G)[{1\over p}]$. Thus
$\dim_{\dbF_p}(\Lie(\calI_{\dbF_p}))\Le\dim(\calI_{B(k)})=\dim(\calI_{\dbQ_p})$.
As $\dim(\calI_{\dbF_p})=\dim(\calI_{\dbQ_p})$, we get that
$\dim_{\dbF_p}(\Lie(\calI_{\dbF_p}))\Le\dim(\calI_{\dbF_p})$. This
implies that
$\dim_{\dbF_p}(\Lie(\calI_{\dbF_p}))=\dim(\calI_{\dbF_p})$ and
that the group $\calI_{\dbF_p}$ is smooth. Thus $\calI$ is a
smooth group scheme over $\Spec(\dbZ_p)$.\endproof

\subsection{The group scheme $\calI_0$.} \label{group} Let $\calI_1$ be the smoothening of $\calI$ defined and proved to exist in [BLR, Ch. 7, pp. 174--175]. We recall that $\calI_1$ is a smooth group scheme of finite type over $\Spec(\dbZ_p)$ equipped with a homomorphism $\calI_1\to\calI$ which is uniquely determined by the following universal property (see [BLR, Ch. 7, 7.1, Thm. 5]):

\medskip
{\bf (i)} if $Y$ is a smooth scheme over $\Spec(\dbZ_p)$, then each morphism $Y\to \calI$ factors uniquely through $\calI_1$.
\medskip

The scheme $\calI_1$ is obtained from $\calI$ through a sequence
of dilatations centered on special fibres (see the paragraph
before [BLR, Ch. 7, 7.1, Thm. 5]) and thus it is an affine scheme
over $\calI$ (cf. the very definition of dilatations; see the
first paragraph of [BLR, Ch. 3, 3.2]). Thus $\calI_1$ is an affine
group scheme over $\Spec(\dbZ_p)$. If $p>2$, then from (i) and
Lemma \ref{p3lemma} we easily get that the homomorphism
$\calI_1\to\calI$ is an isomorphism; thus $\calI_1=\calI$. Let
$\calI_0$ be the unique open subgroup scheme of $\calI_1$ whose
special fibre is the identity component of $\calI_{1\dbF_p}$. Thus
there exists a homomorphism $\calI_0\to\calI$ whose generic fibre
is an isomorphism and moreover we have:

\medskip
{\bf (ii)} the special fibre $\calI_{0\dbF_p}$ of $\calI_0$ is a
smooth, connected, affine scheme.

\subsection{Invariants.} Let $n_{\lambda}\in\dbN$ be the smallest number such
that $(H,\lambda)$ is uniquely determined up to isomorphism by
$(H[p^{n_{\lambda}}], \lambda[p^{n_{\lambda}}])$. Its existence is
implied by [Va1, Subsection 3.2.5]. Let $t_{\lambda}\in\dbN$ be the
$i$-number of $(M,\phi,G)$ defined in [Va1, 3.1.4] (i.e., the
smallest natural number such that for each element $g\in G(W(k))$
congruent mod $p^{t_{\lambda}}$ to $1_M$, there exists an
isomorphism between $(M,g\phi)$ and $(M,\phi)$ which is an element
of $G(W(k))$). From an argument entirely analogous to the proof of
Theorem \ref{supersingularinvariants} (a) (cf. [Va1, Subsections
3.2.1 and 3.2.5]), we get that $n_{\lambda}=t_{\lambda}$.

\subsection{Proof of Theorem \ref{traversopp}.} We will use the notations of Subsection 2.1 to prove that $t_{\lambda}\Le m$. Let $g\in G(W(k))$ be congruent mod $p^m$ to $1_M$. As $g\in\calA(W(k))$ (see proof of Theorem
\ref{supersingularinvariants} (a)) and $\iota(g)g=1_M$, we have
$g\in\calI(W(k))$. We show that in fact we have $g\in\calI_0(W(k))$.

We first show that $g\in\calI_1(W(k))$. If $p>2$, this is obvious
as $\calI_1=\calI$. Suppose that $p=2$. Let $R$ be a
$\dbZ_2$-subalgebra of $W(k)$ of finite type such that the
morphism $\Spec(W(k))\to \calI$ defined by $g$, factors through
$\Spec(R)$. The monomorphism $\dbZ_2\hookrightarrow W(k)$ is of
index of ramification $1$ and the generic point of $\Spec(R)$
belongs to the smooth locus of $\Spec(R[{1\over 2}])$ over
$\Spec(\dbQ_2)$. Based on these and [BLR, Ch. 3, 3.6, Prop. 4], we
get that there exists an $R$-algebra $R_1$ which is smooth over
$\dbZ_2$ and for which there exists an $R$-homomorphism $R_1\to
W(k)$ (in fact, the proof of loc. cit. shows that one can assume
that we have $R_1[{1\over 2}]=R[{1\over 2}]$ and thus that $R_1$
is an $R$-subalgebra of $W(k)$). Thus we can view $g$ as an
$R_1$-valued point of $\calI$. From this and Subsection
\ref{group} (i) we get that we can view $g$ as an $R_1$-valued
point of $\calI_1$. Thus $g\in \calI_1(W(k))$ even if $p=2$.

As for each prime $p$ the number of connected components of $\calI_{1\dbF_p}$ is finite (i.e., the group $\calI_1(k)/\calI_0(k)$ is finite), there exists $s\in\dbN$ such that the images of the two groups
$\Ker(G(W(k))\to G(W_m(k)))$ and $\Ker(G(W_{m+s}(k))\to
G(W_m(k)))$ in $\calI_1(k)/\calI_0(k)$ are
equal. But $\Ker(G(W_{m+s}(k))\to G(W_m(k)))$ is the group of
$k$-valued points of a connected group over $k$ which has a composition series whose factors are isomorphic to the vector group
over $k$ defined by the Lie algebra $\Lie(G_k)$. From the last two sentences we get that the image of $\Ker(G(W(k))\to G(W_m(k)))$ in
the finite group $\calI_1(k)/\calI_0(k)$ is the identity. Thus we have $g\in\calI_0(W(k))$.

As $\calI_0$ is a smooth group scheme over $\Spec(\dbZ_p)$ whose
special fibre is a connected, affine scheme (cf. Subsection
\ref{group} (ii)) and as $g\in\calI_0(W(k))$, from Proposition
\ref{wellknown} we get that there exists an element
$g_0\in\calI_0(W(k))\leqslant G(W(k))$ such that
$g_0^{-1}\sigma(g_0)=g$. Thus $g_0g\sigma(g_0)^{-1}=1_M$. As
$\sigma(g_0)=\phi(g_0)$, we have $g_0g\phi(g_0)^{-1}=1_M$. Thus
$g_0g\phi g_0^{-1}=\phi$ i.e., $g_0\in G(W(k))$ is an isomorphism
between $(M,g\phi)$ and $(M,\phi)$. This implies that
$t_{\lambda}\le m$.

As $n_{\lambda}=t_{\lambda}\Le m$, from Corollary \ref{finaltouch}
we get $n_{\lambda}\Le q+1\Le d$. The inequality
$n_{\lambda}\Le d$ ends the proof of Theorem \ref{traversopp}.
For $p>2$, the inequality $t_{\lambda}\Le m$ refines the
inequality $t_{\lambda}\Le m+1$ which is a particular case of [Va1,
Example 3.3.6].\endproof

\medskip\noindent
{\bf Acknowledgments.}

\medskip
The first author has been supported by the Japanese Society for
the Promotion of Science (JSPS) while working on this paper at the
University of T\=oky\=o. The second author would like to thank
University of Arizona for good conditions with which to write this
note. The authors would like to thank the referee for pointing out
a minor inaccuracy.

%}} \noindent

\bigskip
\hbox{Marc-Hubert Nicole,\;\;\;Email: nicole@ms.u-tokyo.ac.jp}
\hbox{Address: University of T\=oky\=o, Department of Mathematical
Sciences,} \hbox{Komaba, 153-8914, T\=oky\=o, Japan.}

\hbox{}

\hbox{Adrian Vasiu,\;\;\;Email: adrian@math.binghamton.edu; fax:
1-607-777-2450.} \hbox{Address: Department of Mathematical
Sciences, Binghamton University, } \hbox{Binghamton, New York,
13902-6000, U.S.A.}

\end{document}